\documentclass[12pt]{article}
\usepackage{amsfonts}
\usepackage{latexsym}
\usepackage{bm}
\usepackage{psfrag}
\usepackage[dvips]{graphicx}
\usepackage[fleqn]{amsmath}
\usepackage{theorem}

\newtheorem{theorem}{Theorem}
{\theorembodyfont{\rmfamily} \newtheorem{exam}[theorem]{Example}}
\newtheorem{lemma}[theorem]{Lemma}
\newtheorem{lem}[theorem]{Lemma}

\newtheorem{conj}[theorem]{Conjecture}

\def\b{\backslash}

\def\tr{\textrm}

\def\A{\mathcal{A}}
\def\B{\mathcal{B}}

\def\F{\mathcal{F}}
\def\G{\mathcal{G}}
\def\I{\mathcal{I}}
\def\J{\mathcal{J}}

\def\mc{\mathcal}
\newcommand{\m}[1]{\ensuremath{\left\vert #1 \right\vert}}
\newcommand{\Z}{{\ensuremath{\mathbb{Z}}}}  

\begin{document}
\title{Graphs with the Erd\H os-Ko-Rado Property}
\author{Fred Holroyd$^{\!1}$ and John Talbot$^{\!2}$\\[5mm]
$^{1}$Department of Pure Mathematics, The Open University,\\
Walton Hall, Milton Keynes MK7~6AA, United Kingdom\\
\texttt{f.c.holroyd@open.ac.uk}\\ \\
$^{2}$Merton College and the Mathematical Insititute\\
Oxford University, OX1 4JD\\
\texttt{talbot@maths.ox.ac.uk}}
\date{\today} \maketitle
\begin{abstract}
For a graph $G$ and integer $r \ge 1$ we denote the family of
independent $r$-sets of $V(G)$ by $\I^{(r)}(G)$. A graph $G$ is
said to be \emph{$r$-EKR} if no intersecting subfamily of
$\I^{(r)}(G)$ is larger than the largest such family all of whose
members contain some fixed $v \in V(G)$. If this inequality is
always strict, then $G$ is said to be \emph{strictly $r$-EKR}.  We
show that if a graph $G$ is $r$-EKR then its lexicographic product
with any complete graph is $r$-EKR.

For any graph $G$, we define $\mu(G)$ to be the minimum size of a
maximal independent vertex set.  We conjecture that, if $1 \le r
\le \frac{1}{2} \mu(G)$, then $G$ is $r$-EKR, and if $r<
\frac{1}{2}\mu(G)$, then $G$ is strictly $r$-EKR. This is known to
be true when $G$ is an empty graph, a cycle, a path or the
disjoint union of complete graphs. We show that it is also true
when $G$ is the disjoint union of a pair of complete multipartite
graphs.
\end{abstract}
\section{Introduction}
An \emph{independent} set in a graph $G=(V,E)$, is a subset of the
vertices not containing any edges. For an integer $r\geq 1$ we
denote the collection of independent $r$-sets of $G$ by\[
\I^{(r)}(G)=\{A\subset V(G):\m{A}=r\tr{ and $A$ is an independent
set}\}.\] A subfamily $\A$ of $\I^{(r)}(G)$ is said to be
\emph{intersecting} if $A,B \in \A$ implies $A\cap B \ne
\emptyset$. If $v\in V(G)$ then the collection of independent
$r$-sets containing $v$ is \[ \I_v^{(r)}(G)=\{A\in
\I^{(r)}(G):v\in A\}.\] Such a family is called a \emph{star}.

A graph $G$ is said to be \emph{$r$-EKR} if no intersecting family
$\A \subseteq \I^{(r)}(G)$ is larger than the largest star in
$\I^{(r)}(G)$. If every intersecting family
$\A\subseteq\I^{(r)}(G)$ of maximum size is a star then $G$ is
said to be \emph{strictly $r$-EKR}.

The classical result in this area is the Erd\H os-Ko-Rado theorem
which can be stated as follows.

\begin{theorem}[Erd\H os-Ko-Rado \cite{E}]
If $G=E_n$ is the empty graph of order $n$, then $G$ is $r$-EKR if
$n\geq 2r$ and strictly $r$-EKR if $n>2r$.
\end{theorem}

There are several other recent results of this type.

\begin{theorem}[Bollob\'as and Leader\cite{BL}]\label{BLthm}
If $n \geq r$, $t\geq 2$ and $G$ is the disjoint union of $n$
copies of $K_t$, then $G$ is $r$-EKR and strictly so unless $t=2$
and $n=r$.
\end{theorem}

\begin{theorem}[Holroyd and Talbot \cite{HT}]\label{completethm}If $G$ is the disjoint union of $n\geq r$ complete graphs each of order at least two, then $G$ is $r$-EKR.
\end{theorem}

In this paper we consider the question of when a graph is $r$-EKR.
In the next section we give the first of our two main results: if
a graph $G$ is $r$-EKR then its lexicographic product with any
complete graph is also $r$-EKR.

In section 3 we present some examples showing that graphs exhibit
a variety of EKR properties. These serve to motivate a conjecture
we propose, giving a lower bound on the minimum $r$ such that a
given graph $G$ can fail to be $r$-EKR. This conjecture is known
to be true for empty graphs, cycles, paths and disjoint unions of
complete graphs. In the final section we give our second main
result: this conjecture is true for disjoint unions of two
complete multipartite graphs.

Throughout $G$ is assumed to be a simple graph (without loops or
multiple edges) and to have finite vertex set $V(G)$ and edge set
$E(G)$. The independence number of a graph is denoted by
$\alpha(G)$ and the \emph{minimax independence number} (the
minumum size of a maximal independent vertex set) by $\mu(G)$.

An \emph{anomalous} subfamily of $\I^{(r)}(G)$ is an intersecting
subfamily that is not a subfamily of any star. A vertex $v$ is an
$r$-\emph{centre} of $G$ if $\m{\A} \leq \m{\I_v^{(r)}(G)}$ for
every intersecting subfamily $\A$ of $\I^{(r)}(G)$ and is a
\emph{strict} $r$-\emph{centre} if $\m{\A} < \m{\I_v^{(r)}(G)}$
for every anomalous subfamily $\A$ of $\I^{(r)}(G)$.

Where no confusion is caused, we may omit the argument `$(G)$'.


If $\F$ is a family of sets then we define
\begin{align*}\F_x &=\{A\in\F\colon x \in A\},\\
\F^{(r)}&=\{A \in \F\colon \m{A}=r\},\\
\F_x^{(r)}&=\F_x \cap \F^{(r)}.
\end{align*}

Given two graphs $G$ and $H$, the \emph{lexicographic product}
$G[H]$ is constructed (informally speaking) by replacing each
vertex of $G$ with a copy of $H$.  More formally, $V(G[H]) =
V(G)\times V(H)$, where $(v,w)$ is adjacent in $G[H]$ to $(x,y)$
if and only if \emph{either} $v$ is adjacent to $x$ in $G$
\emph{or} $v=x$ and $w$ is adjacent to $y$ in $H$.

It is useful to develop a generalization of this concept: rather
than insisting that each vertex of $G$ be replaced by a copy of a
fixed graph, we may allow the replacement graphs to vary.  For
example, if we begin with $G$ and replace each vertex $v_1,
\ldots, v_k$ with a copy of a graph $H$ and each vertex $w_1,
\ldots, w_q$ with a copy of a graph $J$, then we denote the result
by\\ $G[v_1, \ldots, v_k: H;\mbox{ }w_1, \ldots, w_q: J]$.
\section{Lexicographic products with complete graphs}
We begin with a lemma concerning EKR properties of general set
families, inspired by the elegant proof due to Katona \cite{K} of
the Erd\H os-Ko-Rado Theorem and giving it a more general context.

A family of subsets of a set $S$ is a $q$-\emph{covering} of $S$
if each element of $S$ belongs to exactly $q$ sets of the family.

\begin{lemma}\label{Genkat}
Let $\F$ be a family of $r$-subsets of a finite set $S$, let
$\Gamma$ be a family of subfamilies of $\F$, let $x \in S$, and
suppose suppose that, for some $q$:
\begin{itemize} \item[(i)] $\Gamma$ is a $q$-covering of $\F$;
\item[(ii)]  $x$ is an $r$-centre of each $\G \in \Gamma$.
\end{itemize}
Then $x$ is an $r$-centre of $\F$.
\end{lemma}
\textbf{Proof.\/} Let $\A$ be any intersecting subfamily of $\F$.
Since $\Gamma$ is a $q$-covering of $\F$, it is a $q$-covering of
$\A$ and so
\begin{equation}\label{qF}
q\m{\A} = \sum_{\G\in\Gamma}\m{\A\cap\G}.
\end{equation}
In particular,
\begin{equation}\label{qX}
q\m{\F_x} = \sum_{\G\in\Gamma}\m{\G_x}.
\end{equation}
But for any intersecting subfamily $\A$ of $\F$ and any $\G \in
\Gamma$, the family $\A\cap\G$ is an intersecting subfamily of
$\G$, and so (since $x$ is an $r$-centre of each $\G$)
\begin{equation}\label{XG}
\m{\A\cap\G} \leq \m{\G_x}\mbox{ }(\G\in\Gamma).
\end{equation}
Now, (\ref{qF}), (\ref{qX}) and (\ref{XG}) imply (for any
intersecting subfamily $\A$ of $\F$):
\[\m{\A} \leq \m{\F_x},\]
and so $x$ is an $r$-centre of $\F$. \hfill{$\Box$}

\textbf{Remark.\/} The `strict' extension of Lemma~\ref{Genkat} is
false.  For example, let $S$ be the vertex set of an octahedron
and let $\F$ be the family of 3-subsets of $S$ corresponding to
the faces.  Let $\Gamma$ be the 1-covering (i.e. partition) of
$\F$ into pairs of opposite faces.  Each $\G \in \Gamma$ is
trivially EKR, and so each $x\in S$ is a strict 3-centre of each
such $\G$. Also, each $x \in S$ is a 3-centre of $\F$ with
$\m{\F_x}=4$.  However, there exist anomalous subfamilies of $\F$
of cardinality 4, namely (for each face $F$) the family of faces
containing at least two of the vertices of $F$.  Thus the elements
of $S$ are not strict 3-centres of $\F$. \hfill{$\Box$}

\begin{lem}\label{lexiprod} Let $v$ be an $r$-centre of a graph
$G$ and let $m \in \Z^+$; then each vertex $(v,x)\mbox{ }(x\in
V(K_m)$ is an $r$-centre of the lexicographic product $G[K_m]$.
\end{lem} \textbf{Proof.\/} When $m=1$ the statement is trivial,
so assume $m>1$.

For the purposes of this proof, it is convenient to identify the
vertices (in some fixed way) with the elements of the set
$[n]=1,2,\ldots,n$, and to identify the vertices of $K_m$ with the
elements of the cyclic group $\Z_m$. Let $\F$ be the family of
functions $f\colon [n] \to Z_m$. Then, for each $X \in
\I^{(r)}(G)$ and each $f \in \F$, we define
\[X \circ f=\{(v,f(v))\colon\mbox{ }v \in X\}.\]
We now define an equivalence relation $\sim$ on $\F$ by \[f \sim
g\mbox{ whenever }g=g+z\mbox{ for some }z \in Z_m;\] that is,
$f(v)=g(v)+z\mbox{ }(v \in [n]).$  We denote by $\Psi$ the family
of equivalence classes, and for each $\psi \in \Psi$ we let
$\J_\psi$ denote the following subfamily of $\I^{(r)}(G[K_m])$:
\[\J_\psi = \{X \circ f:\mbox{ }X \in \I^{(r)}(G),\mbox{ }f \in
\psi\}.\] Each $Y \in \I^{(r)}(G[K_m])$ is of the form $X \circ f$
for exactly one $X \in \I^{(r)}(G)$ and exactly $m^{n-r}$
functions $f$ (each in a distinct equivalence class).  That is,
the family $\{\J_\psi:\mbox{ }\psi \in \Psi\}$ is a $q$-covering
of $\I^{(r)}(G[K_m])$ where $q=m^{n-r}$.  By Lemma~\ref{Genkat},
it remains to show that each $(v,x)\mbox{ }(x \in \Z_m)$ is an
$r$-centre of $\J_\psi$ for each $\psi \in \Psi$.

Let $\psi \in \Psi$ and let $\A$ be an intersecting subfamily of
$\J_\psi$.  Let \[\B=\{X \in \I^{(r)}(G):\mbox{ }X \circ f \in
\A\mbox{ for some }f \in \psi\}.\] Then $\B$ is an intersecting
subfamily of $\I^{(r)}(G)$, and so $\m{\B} \leq
\m{\I_v^{(r)}(G)}$. If $X \in \I^{(r)}(G)$ and $f$, $g$ are
distinct elements of $\psi$, then $(X \circ f) \cap (X \circ
g)=\emptyset$.  But $\A$ is intersecting; thus any two distinct
elements of $\A$ correspond to distinct elements of $\B$. Hence
$\m{\A}=\m{\B}$, and so
\begin{equation}\label{lexic}\m{\A} \leq
\m{\I_v^{(r)}(G)}.\end{equation} Let $x \in \Z_m$ and consider the
vertex $(v,x)$ of $G[K_m]$. For each $\psi \in \Psi$ and each $X
\in \I^{(r)}(G)$, we have $(v,x) \in X \circ f\mbox{ for some }f
\in \psi$ if and only if $X \in \I_v^{(r)}(G)$, in which case
there is exactly one $f \in \psi$ with this property.  Thus
$\m{(\J_\psi)_{(v,x)}}=\m{\I_v^{(r)}(G)}$, and it follows from
(\ref{lexic}) that $(v,x)$ is an $r$-centre of $\J_\psi$.
\hfill{$\Box$}

\begin{theorem}\label{lexithm}If $G$ is $r$-EKR and $m\geq 1$ then $G[K_m]$ is $r$-EKR.
\end{theorem}
\textbf{Proof.} This follows directly from Lemma \ref{lexiprod}.
$\hfill \Box$

It is natural to ask whether Theorem~\ref{lexithm} extends to
lexicographic products that involve replacing the vertices of $G$
with complete graphs of variable rather than constant order.  We
now show that this is not unconditionally true.

\begin{exam}\label{thirteen}
Let $G$ be the graph with vertex set $\{v_1,\ldots,v_{13}\}$
depicted below (Figure~\ref{thirteenfig}).
\begin{figure}\label{thirteenfig}
\begin{center}
\psfrag{a}[][]{$v_8$} \psfrag{b}[][]{$v_7$} \psfrag{c}[][]{$v_1$}
\psfrag{d}[][]{$v_4$} \psfrag{e}[][]{$v_{10}$}
\psfrag{f}[][]{$v_{11}$} \psfrag{g}[][]{$v_{12}$}
\psfrag{h}[][]{$v_6$} \psfrag{i}[][]{$v_5$} \psfrag{j}[][]{$v_2$}
\psfrag{k}[][]{$v_3$} \psfrag{l}[][]{$v_9$}
\psfrag{m}[][]{$v_{13}$}
\includegraphics[height=210pt]{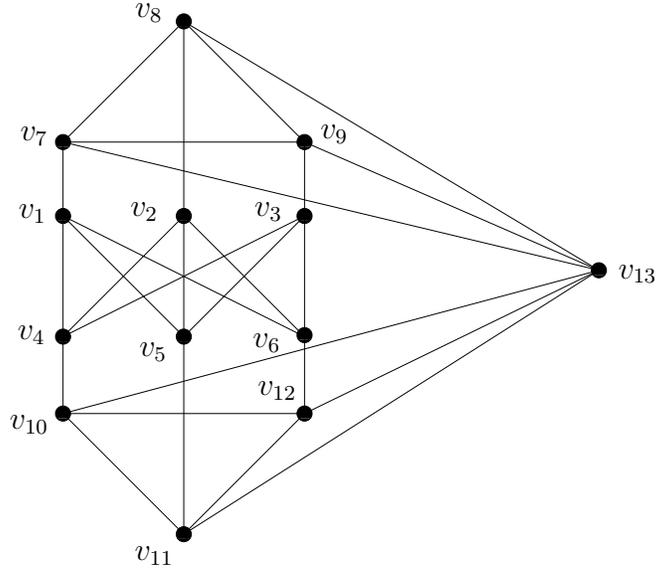}
\caption{the graph $G$ of Example~\ref{thirteen}}
\end{center}
\end{figure}
It may straightforwardly be verified that $G$ is 3-EKR, the
vertices $v_1,\ldots,v_6$ being 3-centres, with
$\m{\I_v^{(3)}(G)}=17\quad(v=v_1,\ldots,v_6).$  The family of
independent vertex 3-sets containing at least two of the vertices
$v_1$, $v_2$, $v_3$ is of cardinality 16 and is one of two
anomalous families of maximum cardinality.

Now let $m \in \Z^+$ and consider the graph $G[v_{13}\colon K_m]$.

Then, $\m{\I_v^{(3)}(G[v_{13}\colon K_m])}=15+2m\mbox{
}(v=v_1,\ldots,v_6)$, the values for the remaining vertices being
independent of $m$.  However, the anomalous family consisting of
all independent 3-sets of $G[v_{13}\colon K_m]$ containing at
least two of the vertices $v_1$, $v_2$, $v_3$ is of cardinality
$13+3m$.  Thus, for $m>2$, the vertices $v_1,\ldots,v_6$ of
$G[v_{13}\colon K_m]$ are not 3-centres (and $G[v_{13}\colon K_m]$
is not 3-EKR).
\end{exam}
\section{Examples of EKR behaviour and a conjecture}
Trivially, any graph is 1-EKR.  The question of when a
(non-complete) graph is 2-EKR is easy to deal with:
\begin{theorem}\label{2-EKR}
Let $G$ be any non-complete graph of order $n$ and with minimum
degree $\delta$.
\begin{itemize}
\item[(i)] If $\alpha=2$, then $G$ is strictly 2-EKR. \item[(ii)]
If $\alpha \geq 3$, then $G$ is 2-EKR if and only if $\delta \leq
n-4$ and strictly so if and only if $\delta \leq n-5$, the
2-centres being the vertices of minimum degree.
\end{itemize}\end{theorem}

\textbf{Proof.\/}  Let $\A$ be an anomalous family of independent
vertex 2-sets.  Then $\m{\A} \geq 3$, and $\A$ must contain the
three 2-subsets of some independent 3-set; but then no other 2-set
can intersect all three of these, and so $\A$ must consist exactly
of the three 2-subsets of an independent 3-set.  Thus:
\begin{itemize}
\item[(i)]  If $\alpha = 2$, then there is no anomalous family of
independent vertex 2-sets, so $G$ is strictly 2-EKR; \item[(ii)]
Otherwise, the anomalous families of independent vertex 2-sets are
all of cardinality 3 and the result follows from the fact that,
for any vertex $v$, \[\m{\I_v^{(2)}}=n-1-d(v).\] \hfill{$\Box$}
\end{itemize}

For the remainder of the paper, then, we concentrate on on the
question: for $3 \le r \le \alpha(G)$, when is $G$ $r$-EKR?

All of the graphs studied in \cite{HT}, including those arising
from reinterpreting \cite{E} and \cite{BL}, are $\alpha$-EKR and
also $\lfloor \alpha/2 \rfloor$-EKR, giving rise to the question:
is this always true? The answer is no, as the following examples
show.

\begin{exam}\label{dodeca}
Let $G$ be the graph of the regular dodecahedron (that is, the
graph whose vertices and edges are those of the dodecahedron).

Then $\alpha = 8$, where $\I^8$ consists of the vertex sets of the
five inscribed cubes of the dodecahedron.  Any pair of these sets
intersects on two (opposite) vertices, but any given vertex
belongs to just two of them.  Thus $\I^8$ is an anomalous family
and $G$ is not 8-EKR. We note, without proof, that if $G$ is the
graph of any of the Platonic solids other than the dodecahedron,
then $G$ $\alpha$-EKR.
\end{exam}
\begin{exam}\label{spiky}
Let $F$ be the graph with vertices $v_1,\ldots,v_7$ where
$v_1,\ldots,v_4$ are pairwise adjacent and $v_{i+4}$ is adjacent
only to $v_i\mbox{ }(i=1,2,3)$.  (See Figure~\ref{spikyfig}.)

Now let $G=F[v_1,v_2,v_3:K_3;\mbox{ }v_4:E_4]$.  Then the order of
$G$ is 16 and $\alpha=7,\mbox{  }\mu=3$.  Moreover, the families
$\I^{(r)}\mbox{ }(4\leq r\leq 7)$ are precisely the families of
$r$-subsets of the unique independent 7-set.  Thus, $G$ is 7-EKR
in a trivial way and (by the Erd\H os-Ko-Rado Theorem) is not 4-,
5- or 6-EKR.  More interestingly, $G$ fails to be 3-EKR, since no
vertex belongs to more than 21 independent 3-sets but there is an
anomalous family consisting of the 22 independent 3-sets
containing at least two of $v_5$, $v_6$, $v_7$. Thus it is
possible for a graph to fail to be $\lfloor \alpha/2 \rfloor$-EKR
and to fail to be $\mu$-EKR.
\begin{figure}\label{spikyfig}
\begin{center}
\psfrag{a}[][]{$v_{5}$} \psfrag{b}[][]{$v_1$}
\psfrag{c}[][]{$v_4$} \psfrag{d}[][]{$v_2$} \psfrag{e}[][]{$v_6$}
\psfrag{f}[][]{$v_3$} \psfrag{g}[][]{$v_{7}$}
\includegraphics[height=150pt]{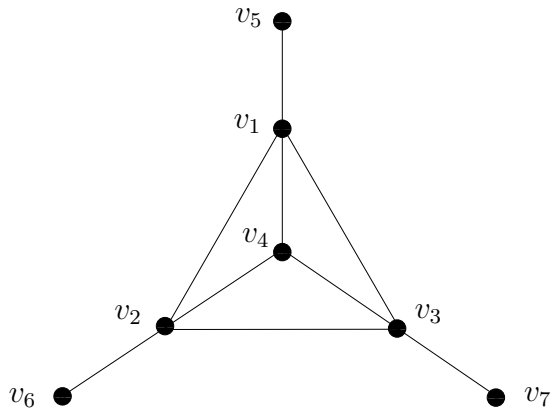}
\caption{the graph $F$ of Example~\ref{spiky}} \label{spikyfig}
\end{center}
\end{figure}
\end{exam}

In each graph studied so far, when $G$ is $\alpha$-EKR, it is so
in a trivial way; but this is not so in general, as the next
example shows.

\begin{exam}\label{icosa}
Let $G$ be the graph of the regular icosahedron. Then $\alpha =
3$. It is straightforward to check that $\m{\I_v^{(3)}}=5$ for any
vertex $v$, and with a little care it is possible to construct an
anomalous family of four independent 3-sets and to verify that no
such family can be extended to a fifth member.  Thus $G$ is
(strictly) 3-EKR.

Note that the antipodal pairs of vertices of $G$ are maximal
independent sets, so that $\mu = 2$.  Therefore, this example also
shows that it is possible for a graph to be $r$-EKR for some
$r>\mu$.
\end{exam}

It is easy to vary Example~\ref{spiky} to produce a graph of
arbitrarily large independence number that fails to be 3-EKR
since, if we replace $K_3$ by $K_p$ and $E_4$ by $E_q$ in the
generalized lexicographic construction of that example, then
$\alpha=q+3$ , the maximum value of $\m{\I_v^{(3)}}$ is
\[\max\{1+2(p+q)+\frac{1}{2}q(q-1),\mbox{
}\frac{1}{2}(q+1)(q+2)\},\] and there is an anomalous subfamily of
$\I^{(3)}$ of cardinality $1+3(p+q)$.  More generally it is
possible, for any $r\geq 3$, to produce a graph of arbitrarily
large independence number that fails to be $r$-EKR.  However, this
does not seem to be true for the minimax independence number.  We
make the following conjecture.
\begin{conj}\label{mu}
Let $G$ be any graph and let $1 \leq r \leq \frac{1}{2}\mu$; then
$G$ is $r$-EKR (and is strictly so if $2 < r < \frac{1}{2}\mu$).
\end{conj}

Each of the above bounds is sharp, as our final example shows.

\begin{exam}\label{bipexam}
Let $G$ be the disjoint union of two copies of the complete
bipartite graph $K_{3,3}$.  Then (by Theorem~\ref{bip} of
Section~4) $\mu=6$ and $G$ is non-strictly 3-EKR and strictly
2-EKR, but not 4-EKR.
\end{exam}
\section{Unions of complete multipartite graphs}
It seems plausible that if any graphs fail to be $r$-EKR, for some
$r \le \frac{1}{2}\mu$, then the smallest examples should have
$\mu = \alpha$ (that is, all maximal independent vertex sets
should have the same cardinality).  This motivates the study of
classes of graphs with this property.

Conjecture~\ref{mu} is already known to hold for certain classes
of graphs; in particular it holds for empty graphs and disjoint
unions of complete graphs (both of which have $\mu=\alpha$). We
now show that the conjecture also holds for the class of unions of
pairs of complete multipartite graphs; moreover, if $G$ is such a
graph and $\mu=\alpha$, then the bound is sharp in that $G$ fails
to be $r$-EKR if $\frac{1}{2}\mu < r < \mu$.

\begin{theorem}\label{bip}
Let $G$ be a union of two complete multipartite graphs; then:
\begin{itemize}
\item[(i)] $G$ is $r$-EKR if $1 \le r \le \frac{1}{2}\mu$; \item[(ii)]$G$ is strictly $r$-EKR if $2 < r <
\frac{1}{2}\mu$;\item[(iii)] $G$ is not $r$-EKR if $\mu =\alpha$
and $\frac{1}{2}\mu < r < \mu$.\end{itemize}
\end{theorem}

Before proving this result, we require further notation and
lemmas.

Let $b_1 \ge b_2 \ge \ldots \ge b_a$. We denote by
$K_a[b_1,b_2,\ldots , b_a]$ the complete $a$-partite graph with
partite sets of sizes $b_1, b_2,\ldots , b_a$ respectively.

Let $G$ be the disjoint union of two complete bipartite graphs,
$G_1=K_a[b_1,\ldots, b_a]$ and $G_2 = K_c[d_1,\ldots d_c]$. Denote
the partite sets of $G_1$ by $V_1,\ldots V_a$ where
$V_i=\{v_{i,1},\ldots,v_{i,b_i}\}\quad(i=1,\ldots,a)$ and those
of $G_2$ by $W_1,\ldots,W_c$ where
$W_i=\{w_{i,1},\ldots,w_{i,d_i}\}\quad(i=1,\ldots,c)$.

For $2 \le i \le a$, define $\phi_i\colon V(G)\to V(G)$ as
follows.\begin{align*}\phi_i(v_{i,j}) &=v_{1,j}\quad&(v_{i,j}\in V_i),\\
\phi_i(v) &=v\quad&\mbox{(otherwise).}\end{align*} Similarly, for
$2 \le i \le c$, define $\theta_i\colon V(G)\to V(G)$ by
\begin{align*} \theta_i(w_{i,j}) &=w_{1,j}\quad&(w_{i,j}\in W_i),\\
\theta_i(w) &=w\quad&\mbox{(otherwise).}\end{align*} With slight
abuse of notation, if $A \in \I(G)$, we may write\\
$\phi_i(A)=\{\phi_i(x)\colon x \in A\}$ and
$\theta_i(A)=\{\theta_i(x)\colon x \in A\}$.  Note that\\
$\phi_i(A), \theta_i(A) \in \I(G)$. We now define the
\emph{compressions} $\Phi_i, \Theta_i$ on subfamilies of $\I(G)$
as follows.  Let $\A \subseteq \I(G)$ and let $2 \le i \le a$.
Then
\[\Phi_i(\A)=\{\theta_i(A):A\in\A\}\cup\{A:A,\theta_i(A)\in\A\}.\]
More informally, for each $A \in \A$ that intersects $V_i$, we
replace $A$ by $\phi_i(A)$ provided that $\phi_i(A)$ is not
already in $\A$; otherwise, we leave $A$ alone.

The compressions $\Theta_i\quad(2 \le i \le c)$ are similarly
defined.

We now note that, if $\A$ is a non-empty intersecting subfamily of
$\I^(G)$, then there is some partite set of $G_1$ or $G_2$ that
intersects every set of $\A$; for any $A \in \A$ is a subset of
$V_i \cap W_j$ for some $i,j$ and now there cannot be $B,C\in \A$
with $B\cap V_i=\emptyset$ and $C\cap W_j=\emptyset$. By
exchanging $G_1$ and $G_2$ if necessary, we may assume that some
fixed $V_i$ intersects each set of $\A$. Clearly, $\B=\Phi_i(\A)$
is an intersecting family with $\m{\B}=\m{\A}$ such that $V_1$
intersects each set of $\B$. Thus, in investigating the sizes of
intersecting subfamilies $\A$ of $\I^{(r)}(G)$, we may assume that
$V_1$ intersects each $A\in \A$; such a family is said to be
\emph{standardized}.

Our first lemma says that any compression of a standardized
intersecting family in $\I^{(r)}(G)$ is a standardized
intersecting family of the same size.

\begin{lemma}\label{firstlem}
Let $2 \le i \le c$.  With the above notation, if $\A\subseteq
\I^(G)$ is standardized and intersecting then so is
$\Theta_i(\A)$, and $\m{\Theta_{i}(\A)}=\m{\A}$.
\end{lemma}
\textbf{Proof.\/} It follows immediately from the definitions that
$\Theta_{i}(\A)$ is standardized and that
$\m{\Theta_{i}(\A)}=\m{\A}$. We now show that $\Theta_{i}(\A)$ is
intersecting.

Let $A,B\in\Theta_i(\A)$. If $A,B\in\A$ then $A\cap B \neq
\emptyset$. Also if $A=\theta_i(C)$ and $B=\theta_i(D)$, with
$C,D\in\A$ and $A,B\notin \A$, then $C\cap D\neq \emptyset$,
implying that $A\cap B \neq \emptyset$. So we may suppose that
$A\in \A\cap \Theta_i(\A)$ and $B\in\Theta_i(\A)\b\A$.

$A\in\A\cap \Theta_i(\A)$ implies that $C=\theta_i(A)\in\A$. Also
$B\in\Theta_i(\A)\b \A$ implies that there exists $D\in \A$ such
that $B=\theta_i(D)$. Now if $A\cap D\subseteq W_i$ then $C\cap
D=\emptyset$, a contradiction, since $C,D\in\A$. So there exists
$x\in (A\cap D)\b W_i$. But then $x\in A\cap B$ as required. Hence
$\Theta_i(\A)$ is intersecting. \hfill $\Box$

A family $\B\subseteq \I(G)$ is \emph{compressed} if $\B$ is
fixed under every compression.


\begin{lemma}\label{downlem}
Let $G$ be as above. If $\A\subseteq \I(G)$ is a standardized
intersecting family, then there is a standardized compressed
intersecting family $\B\subseteq\I(G) $ such that $\m{\A}=\m{\B}$
and
 $A\cap B\cap (V_1\cup W_1)\neq \emptyset\quad(A,B \in \B)$.
\end{lemma}
\textbf{Proof.\/} Let $\B = \Theta_2 \circ \Theta_3 \circ \ldots
\circ \Theta_c(\A)$.  Then, for any $A \in \B$ such that $A
\subseteq V_1 \cup W_i$ and $i>1$, we have $\theta_i(A) \in \B$,
and so $\B$ is compressed.  By Lemma~\ref{firstlem}, $\B$ is
intersecting and $\m{\B} = \m{\A}$.  Now let $A,B\in \B$. Suppose
$A \cap B \subseteq W_i$ where $i > 1$.  Then $A \cap \theta_i(B)
= \emptyset$, giving
a contradiction since $A, \theta_i(B) \in \B$.  \hfill $\Box$\\

\textbf{Proof of Theorem~\ref{bip}}

\textbf{Proof of (i)} Let $G = G_1 \cup G_2$ as above and let
$\A\subseteq \mc{I}^{(r)}(G)$ be an intersecting family.  By using
Theorem~\ref{2-EKR} or by direct consideration of small cases, we
may assume $r \ge 3$. We shall show that
\[\m{\A} \le \m{\I_x^{(r)}}$ for some $x\in V(G).\]

We may assume that $\A$ is standardized; by Lemma \ref{downlem} we
may also assume that $\A$ is compressed and that $A \cap B
\cap(V_1 \cup W_1) \neq \emptyset\quad(A,B\in \A)$.

Partition $\A$ as $\A=\A_0 \cup \A_1 \cup \ldots \cup \A_c$ where
$\A_0 = \{A\in \A\colon A \subseteq V_1\}$ and, for $1 \le i \le
c$, \[\A_i=\{A \in \A\colon A \cap W_i \ne \emptyset\}.\]
Correspondingly, let $\mc{J}=\I_x^{(r)}$ where $x=v_{1,1}\in V_1$,
and partition $\mc{J}$ as $\mc{J}_0 \cup \mc{J}_1 \cup\ldots \cup
\mc{J}_c$.  Now,
\begin{equation}\label{mu} \mu(G) = b_a+d_c \le \m{V_1\cup W_1}=b_1+d_1.\end{equation}
Thus, by the Erd\H os-Ko-Rado Theorem (since $r \le
\frac{1}{2}\mu$), we have
\begin{equation}\label{A0+A1}
\m{\A_0} + \m{\A_1} \le
{{b_1+d_1-1}\choose{r-1}}=\m{\mc{J}_0}+\m{\mc{J}_1}.\end{equation}
%
%
We now compare $\m{\A_i}$ with $\m{\mc{J}_i}\quad(2 \le i \le c)$.
Since each $A$ in $\A_i\cup \mc{J}_i$ intersects $V_1$ and $W_i$,
we have
\[s_i \le \m{A \cap V_1} \le t\quad(A \in \A_i\cup\mc{J}_i)\] where
$s_i=\max\{1,r-d_i\}\mbox{ , } t=\min\{r-1,b_1\}$.

For $2 \le i \le c\mbox{ , }s_i \le j \le t$, let
$\A_i^{(j)}=\{A\in \A_i\colon \m{A\cap V_1}=j\}$ and
$\B_i^{(j)}=\{A\cap V_1\colon A \in \A_i^{(j)}\}$.\\
Analogously, let $\mc{J}_i^{(j)}=\{A\in \mc{J}_i\colon \m{A\cap
V_1}=j\}$ and $\mc{K}_i^{(j)}=\{A\cap V_1\colon A \in
\mc{J}_i^{(j)}\}$. Then, for $2 \le i \le c$:
\begin{equation}\label{Ai}\m{\A_i} \le
\sum_{j=s_i}^t{\m{\B_i^{(j)}}{d_i\choose{r-j}}},\end{equation}
\begin{equation}\label{Ji}\m{\mc{J}_i}=\sum_{j=s_i}^t{\m{\mc{K}_i^{(j)}}{d_i\choose{r-j}}}=\sum_{j=s_i}^t{{{b_1-1}\choose{j-1}}{d_i\choose{r-j}}}.\end{equation}
Since $\A$ is standardized and compressed, each $\B_i$ is
intersecting, by Lemma~\ref{downlem}. Thus, by (\ref{mu}) and the
Erd\H os-Ko-Rado Theorem, we have for $2 \le i \le c,\mbox{  }s_i
\le j \le \frac{1}{2}b_1$:
\begin{equation}\label{Bij}\m{\B_i^{(j)}} \le {{b_1-1}\choose{j-1}}.\end{equation}
Thus, if $t \le \frac{1}{2}b_1$, then we may conclude from
(\ref{A0+A1})--(\ref{Bij}) that $\m{\A} \le \m{\mc{J}}=
\m{\I_x^{(r)}}$.

Suppose now that $t > \frac{1}{2}b_1$. For $s_i \le j \le
\lfloor\frac{1}{2} b_1 \rfloor,\mbox{  }b_1-j \le t$, we have
\begin{equation}\label{A-0} \m{\A_i^{(j)} \cup \A_i^{(b_1-j)}} \le
\m{\B_i^{(j)}}{d_i\choose{r-j}}+\m{\B_i^{(b_1-j)}}{d_i\choose{r-(b_1-j)}}.\end{equation}
Moreover, by the intersecting property, no set in $\B_i^{(b_1-j)}$
can be the complement of a set in $\B_i^{(j)}$, and hence
\begin{equation}\label{Bs} \m{\B_i^{(j)}} + \m{\B_i^{(b_1-j)}} \le
{b_1\choose j}.\end{equation} Two cases arise.

\emph{Case 1} $\m{\B_i^{(b_1-j)}} \le
{{b_1-1}\choose{b_1-j-1}}={{b_1-1}\choose j}$.  Then,
\begin{equation}\label{A-1} \m{\A_i^{(j)} \cup \A_i^{(b_1-j)}} \le
{{b_1-1}\choose{j-1}}{d_i\choose{r-j}}+{{b_1-1}\choose{j}}{d_i\choose{r-(b_1-j)}}=\m{\mc{K}_i^{(j)}\cup
\mc{K}_i^{(b_1-j)}}.\end{equation} \emph{Case 2}
$\m{\B_i^{(b_1-j)}}
> {{b_1-1}\choose j}$.

Now, from $r \le \frac{1}{2}(b_1 + d_i)$, it is straightforward
to deduce that
\begin{equation}\label{dchooses} {d_i\choose{r-(b_1-j)}} \le
{d_i\choose{r-j}}.\end{equation} Together with
inequality~(\ref{Bs}), this implies that (\ref{A-1}) still holds.
Thus,\\ $\m{\A_i} \le \m{\mc{J}_i}\quad(2 \le i \le c)$.  With
equation~(\ref{A0+A1}), this gives the result $\m{\A} \le
\m{\I_x^{(r)}}$, as required.

\textbf{Proof of (ii)}

We now show that $G$ is strictly $r$-EKR for $r < \frac{1}{2}\mu$.

First note that, for $\m{\A}=\m{\I_x^{(r)}}$, equality must hold
in each of the inequalities (\ref{A0+A1}), (\ref{Ai}),
(\ref{Bij}), (\ref{Bs}); moreover, for $r < \frac{1}{2}(b_1+d_i)$,
the inequality (\ref{dchooses}) is strict.  Thus, by the Erd\H
os-Ko-Rado Theorem, for $\m{\A}=\m{\I_x^{(r)}}$, $\A_0 \cup
\A_1=\{A\subseteq V_1\cup W_1\colon \m{A}=r, x_1\in A\}$ for some
$x_1 \in V_1$.  By (\ref{Ai}), (\ref{Ji}), the Erd\H os-Ko-Rado
Theorem now implies, for $2 \le i \le c$: \[\A_0 \cup \A_i
=\{A\subseteq V_1\cup W_i\colon \m{A}=r, x_i \in A\}\] for some
$x_i\in V_1$.  Clearly, the $x_i$ must all be equal and the result
follows.

\textbf{Proof of (iii)}

Since $\mu(G) = \alpha(G)$, there exist $b, d$ such that
$\mu=b+d$, $b_i=b\quad(1\le i \le a)$ and $d_i=d\quad(1\le i \le
c)$.  Let $\frac{1}{2}(b+d) < r < b+d$.

If $r < b$, let $U \subseteq V_1$ such that $\m{U}=r, x \notin U$.
Now $s + (r-1)\ge 2r-d-1>b-1$, so that $U$ intersects every set of
$\I_x^{(r)}$, which is therefore not a maximal intersecting
family.

If $r \ge b$, then let $x \in V_1$ and consider the family $\mc{J}
= (\I_x^{(r)}\b\{A\in \I^{(r)}\colon A\cap V_1 = \{x\}\})\cup \{A
\in \I^{(r)}\colon A\cap V_1=V_1\b\{x\}\}$. It is straightforward
to check that $\mc{J}$ is anomalous, intersecting and larger than
$\I_x{(r)}$. \hfill $\Box$

\end{document}